\documentclass[11pt]{amsart}
\usepackage{amssymb}

\newtheorem{theorem}{Theorem}[section]
\newtheorem{proposition}[theorem]{Proposition}
\newtheorem{lemma}[theorem]{Lemma}
\newtheorem{corollary}[theorem]{Corollary}

\newcommand{\h}{\hbar}

\theoremstyle{definition}
\newtheorem{definition}[theorem]{Definition}

\theoremstyle{remark}

\numberwithin{equation}{section}

\newcommand{\be}{\begin{equation}}
\newcommand{\ee}{\end{equation}}

\newcommand{\bbR}{{\mathbb R}}
\newcommand{\bbN}{{\mathbb N}}

\newcommand{\calV}{{\mathcal V}}

\newcommand{\calH}{{\mathcal H}}

\newcommand{\calC}{{\mathcal C}}
\newcommand{\calR}{{\mathcal R}}
\newcommand{\calS}{{\mathcal S}}

\newcommand{\calZ}{{\mathcal Z}}

\newcommand{\calZtr}{\calZ^\triangle}
\newcommand{\Zperp}{Z^{\circ}_\rho}
\newcommand{\hphi}{\hat\phi}

\DeclareMathOperator{\End}{End}

\DeclareMathOperator{\Hess}{Hess}

\newcommand{\inner}[2]{\langle#1,#2\rangle}
\newcommand{\bra}[1]{\langle#1|}
\newcommand{\ket}[1]{|#1\rangle}
\newcommand{\norm}[1]{\left\lVert#1\right\rVert}

\newcommand{\isom}{\cong}
\newcommand{\del}{\partial}

\renewcommand{\Box}{\square}

\newcommand{\hffm}{{\textstyle \bigwedge^{1/2}}}

\newcommand{\trk}{\text{Tr}^+ K}

\begin{document}

\title[Semiclassical structure of low-energy states]{The semiclassical structure of 
low-energy states in the presence of a magnetic field}
\author{David Borthwick}
\address{Department of Mathematics and Computer Science\\Emory University\\ 
Atlanta, GA 30322}\email{davidb@math.emory.edu}\author{Alejandro Uribe}
\address{Mathematics Department\\
University of Michigan\\ Ann Arbor, Michigan 48109}
\email{uribe@math.lsa.umich.edu}
\thanks{D.B. supported in part by NSF grant DMS-0204985.}
\thanks{A.U. supported in part by NSF grant DMS-0070690.}

\date{December 10, 2004}

\begin{abstract}
We consider a compact Riemannian manifold with a Hermitian line bundle 
whose curvature is non-degenerate.   The Laplacian acting on high tensor powers 
(the semiclassical regime) of the bundle exhibits a cluster of low-energy states.  
We demonstrate that the orthogonal projectors onto these states
are the Fourier components of an operator with the structure of the Szeg\"o
projector, i.e. a Fourier integral operator of Hermite type.  This result yields
semiclassical asymptotics for the low-energy eigenstates.
\end{abstract}

\maketitle

\section{Introduction}\label{intro}

Let $(X,g)$ be a compact Riemannian manifold and $L\to X$ a Hermitian
line bundle with connection, $\nabla$, which we think of a magnetic potential.
The curvature (or field strength) of $\nabla$ will be denoted by $\omega$,
a closed two-form on $X$.  In this paper we will work under the assumption that
\begin{equation}
\omega\ \mbox{is non-degenerate.}
\end{equation}
This assumption implies that there is a bundle automorphism, $K:TX\to TX$
such that
\begin{equation}\label{Kdef}
\forall x\in X\,\ u,v\in T_xX \qquad g_x(u,K(v)) = \omega_x(u,v).
\end{equation}
It is easy to check that $K$ is skew adjoint with respect to $g$, and therefore 
for each $x\in X$ the 
eigenvalues of $K_x$ can be written in the form $\pm i \kappa_j(x)$ where
 $\kappa_j(x)>0$,
$j=1,\ldots , n$ and the dimension of $X$ is denoted: $2n$.  Introduce the function
\begin{equation}\label{03}
\trk(x) = \sum_{j=1}^n\kappa_j(x).
\end{equation}
It is easy to see that this function is smooth, although the individual 
$\kappa_j$ may cross.  
For future reference we also introduce here the associated almost-complex structure, 
$J = K\circ(K^*K)^{-1/2}$.  

We will be concerned with the eigenstates of low energy of the
sequence of operators
\begin{equation}
\Box_k:= -{\nabla^{(k)}}^{*}\nabla^{(k)} - k\>\trk,
\end{equation}
where $\nabla^{(k)}$ is the connection on the $k$-th tensor power
$L^k\to X$.  
The parameter $1/k=\h$ plays the role of Planck's constant, and
the operator $\frac{1}{k^2}\Box_k$ is a perturbation of a
magnetic Schr\"odinger operator by a suitable potential term
times $\h$.  The precise meaning of ``low-energy states"
is the one implied by the following Lemma:

\begin{lemma}\cite{GU, MM}\label{Drift}
There exists a constants $\epsilon, M>0$ such that for large $k$
the spectrum of $\Box_k$ is contained in 
\[
(-\epsilon,\epsilon)\cup[Mk,\infty).
\]
\end{lemma}
A recent proof of this fact by X. Ma and G. Marinescu, \cite{MM}, 
uses only the Lichnerowicz formula for the Spin-c Dirac operator
associated with the almost-K\"ahler manifold $(X,\omega, J)$.

Let $\calH_k\subset L^2(L^k)$ be the span of the low-lying eigensections of $\Box_k$ 
identified by the Lemma. 
The K\"ahler case corresponds precisely to $K=J$ and $J$ integrable.
If in addition $L$ has a complex structure we can
identify (noting that $\trk = n$) $-{\nabla^{(k)}}^{*}\nabla^{(k)} - nk = (1/4) 
\bar\partial^* \bar\partial$.  Thus $\calH_k$ would consist precisely of 
the holomorphic sections of $L^k$.

It is possible to start with a symplectic manifold $(X,\omega)$, and
a Hermitian line bundle with connection $L\to X$ whose curvature is $\omega$.
Then, by choosing a compatible almost-complex structure, $J$, one is
led to the almost-K\"ahler case ($K=J$ but not integrable) 
and the $\calH_k$ are defined for
$k$ sufficiently large.  The \textit{almost-K\"ahler quantization} scheme we
proposed in \cite{BU1} consists of taking $\calH_k$ as the quantum Hilbert space 
corresponding to $\hbar = k^{-1}$.  This is a direct analog of the standard K\"ahler
quantization.

Our main result is a precise description of the
asymptotics, as $k\to\infty$, of the sequence of orthogonal projectors
\[
\Pi_k: L^2(L^k)\to \calH_k.
\]
We will show that these projectors have the same structure as in the
K\"ahler case.  To state our results
more precisely, we introduce the circle bundle $P\to X$, $P\subset L^*$.
Then the sections of $L^k\to X$ are naturally identified with the space of
functions on $P$ that transform under the circle action on $P$ by 
multiplication by the character $e^{ik\theta}$.  Thus
we can regard each $\calH_k$ as a subspace of $L^2(P)$, and extend
the orthogonal
projection to: $\Pi_k:L^2(P)\to \calH_k$. 
The analog of the Szeg\"o projector from the
K\"ahler case would then be the full projection $\Pi = \sum \Pi_k$.
Our result is that $\Pi$ has the microlocal structure of the Szeg\"o projector.
A central role will be played by the submanifold
$\calZ\subset T^*P\setminus\{ 0\}$ of the punctured cotangent
bundle of $P$ consisting of those covectors that are orthogonal to
the conormal to the fibers of
the natural projection $P\to X$.  ($\calZ$ consists of the vertical
covectors.)   Our main result is the following: 

\begin{theorem}\label{One}
The projector $\Pi$ defined above is an Hermite FIO with Schwarz kernel in the
class $I^{1/2}(P\times P,\calZtr)$, where $\calZtr$ is the diagonal relation 
$\{(\zeta,\zeta); \zeta\in\calZ\}\subset T^*(P\times P)$.
\end{theorem}

Results on semiclassical asymptotics for quantized K\"ahler manifolds
have been derived in \cite{BMS,BPU1, BPU2, BU1,BU2}, using 
the structure of the Szeg\"o projector as an Hermite FIO.
Theorem \ref{One} allows immediate extension of these results to the case of
``almost-K\"ahler'' quantization, i.e. quantization of a symplectic manifold 
with compatible almost structure using the spectral projector $\Pi$ defined above
as an analog of the Szeg\"o projector.   All of these are based on the same
principle:  the singularities of the Schwarz kernel $\Pi$ described in Theorem
\ref{One} correspond to the large $k$ behavior of the (finite-rank) projectors
$\Pi_k$.

To understand this connection, one can take a Fourier decomposition of the kernel of 
$\Pi$ as in Lemma 3.4 of \cite{BU2}.  The result is an asymptotic expansion for 
the kernel of $\Pi_k$ (written w.r.t. some local coordinate system and 
trivialization of $L^k$) of the form:
\begin{equation}\label{piexpansion}
\Pi_k(x,y) \sim \sum_{m\in-\bbN_0/2} k^{n-m} e^{ik\theta(x,y)} f_m(x,y,\sqrt{k}(x-y)),
\end{equation}
where the phase $\theta(x,y)$ is essentially determined by the relation $\calZtr$
and $f_m(x,y,u)$ is rapidly decreasing as a function of $u\in \bbR^{2n}$.  
In the process of proving Theorem \ref{One} we will also compute the symbol of $\Pi$
as a Hermite FIO, which is the invariant object associated to the leading term $f_0$
in the expansion above.  The symbol is essentially identical to that of the Szeg\"o
projector in the K\"ahler case.

If we denote by $\psi_j^{(k)}$, $j=1,\dots,d_k$, the (normalized) low energy 
eigenfunctions of
$\Box_k$ which span $\calH_k$, then the  (\ref{piexpansion}) gives an asymptotic
expansion for the orthogonal projection onto these states.  In particular,
taking the symbol into account, we can extract the following leading behavior
(see Theorem 4.1 of \cite{BU2}):
$$
\sum_{j=1}^{d_k} \psi_j^{(k)}(x) \>\overline{\psi_j^{(k)}(y)} 
\sim \left(\frac{k}{2\pi}\right)^{n} 
e^{ik\theta(x,y)} e^{-\frac{k}2 d(x,y)^2} + O(k^{n-1/2}),
$$
uniformly in $x,y$, where the linearization of the phase $\theta(x,y)$ at 
$y=x$ is the bilinear form associated to $K_x$ and $d(x,y)$ is the 
Riemannian distance. 

In the symplectic context, the existence of projectors on $L^2(P)$ with the 
microlocal structure of the Szeg\"o projector (generalized Toeplitz structures) 
was established in the Appendix of \cite{BG}.  
Theorem \ref{One} allows us to realize these structures in an extremely natural
way, i.e. through spectral projection of the Laplacian.

\medskip
\textit{Plan of the proof.}
The analysis is carried out {\em naturally} on the circle bundle, $P$.
Endow $P$ with the Kaluza-Klein metric and
denote by $\Delta_P$ the Laplace-Beltrami operator on $P$. 
The horizontal Laplacian, $\Delta_H$, is the operator defined 
by the identity
\begin{equation}\label{1h}
\Delta_P = D_\theta^2 + \Delta_H.
\end{equation}
Here $D_\theta = -i\partial/\partial\theta$ is differentiation with respect
to the vector field generating the circle action on $P$.  Since the 
circle action is by isometries, $\Delta_P$ commutes with $D_\theta$ and
hence with $\Delta_H$.  Under the natural identification of the sections
of $L^k$ with the $k$-the eigenspace of $D_\theta$, the Laplacian
$-\nabla^{(k)^*}\nabla^{(k)}$ corresponds to the restriction of
$\Delta_H$ to that subspace.  Therefore, under said correspondence
the operators, $\{\Box_k\}$, ``roll up" to the operator
\[
\Delta_H - (\trk) D_\theta
\]
where we continue to denote by $\trk$ the pull-back of (\ref{03}) to $P$.
(Note that the associated multiplication operator commutes with $D_\theta$.)
We will construct the projector, $\Pi$, by taking suitable functions of
the first order operator
\begin{equation}\label{0b}
(D_\theta)^{-1}\,\Delta_H - \trk
\end{equation}
A complication is that this operator is a singular $\Psi$DO (because $D_\theta$
is not elliptic), but the singularities occur away from the characteristic
variety of $\Delta_H$.  We will therefore in fact work with an operator, $A$,
which is microlocally equal to (\ref{0b}) away from the characteristic 
variety of $D_\theta$.
$A$ is a standard $\Psi$DO with double symplectic characteristics
and such that the spectral projector of $A$
onto an interval $(-a,a)$ for small $a$ differs from $\Pi$ by an operator
of finite rank. 

The construction of the projector $\Pi$  is based on the following result:

\begin{theorem}\label{Main}  
Let $A$ be an operator microlocally equal to
the operator (\ref{0b}) in a conic neighborhood of the
characteristics $\calZ$, and let
$\phi\in C_0^\infty(\bbR)$.  Then the operator
\[
\hphi(A) := \int_\bbR\, e^{-itA}\,\phi(t)\,dt
\]
is an Hermite FIO with Schwarz kernel in the class $I^{1/2}(P\times P,\calZtr)$.  
\end{theorem}
We will identify the symbol 
of $\hphi(A)$ with an object that associates to each point of $\calZ$
the following function of  the harmonic oscillator on the symplectic normal
bundle of $\calZ$:
\[
\sum_{\nu=(\nu_1,\ldots,\nu_n)}\,\hat{\phi}(2\kappa\cdot\nu)\,
\ket{\nu}\bra{\nu}\,.
\]
With a suitable choice of test function, $\phi$, the operator $\hphi(A)$ 
is close to the spectral projector $\Pi$.  A Neumann series argument finishes
the construction of $\Pi$, proving Theorem \ref{One}.  The construction
shows furthermore that the symbol 
of $\Pi$ is projection onto the ground state of the harmonic oscillator on the 
symplectic normal bundle of $\calZ$:
$$
\sigma(\Pi) = \ket{0}\bra{0},
$$
i.e. the the symbol is the same as for the Szeg\"o projector in the K\"ahler case.

As we were finishing writing this paper Xiaonan Ma and George Marinescu  
forwarded us their preprint, \cite{MM2}, where they obtain, among
other things, a description of the projectors $\Pi_k$ (by methods
quite different from ours).

\section{Hermite Fourier integral operators}\label{FIOsec}

An Hermite FIO differs from a standard FIO in having a canonical relation 
that is an isotropic, rather than Lagrangian, submanifold of the cotangent bundle.
 This generalization was motivated by the structure of the Szeg\"o projector
on the boundary of a strictly pseudoconvex domain.

For completeness, we begin by recalling the definition of a Hermite distribution.
Let $M$ be a smooth manifold and $\Sigma \subset T^*M\backslash\{0\}$ 
a homogeneous isotropic submanifold.  

Let $B$ be an open conic subset of $(\bbR\times
\bbR^n)\backslash\{0\}$, given the coordinates $(\tau, 
\eta)$.   A \textit{non-degenerate phase function} is a function 
$\psi\in C^\infty(M\times B, \bbR)$ which satisfies:
\begin{enumerate}
\item $\psi(x,\tau,\eta)$ is homogeneous in $(\tau,\eta)$.
\item $d\psi$ is nowhere zero.
\item The critical set of $\psi$,
$$
C_\psi = \{(x,\tau,\eta);  (d_\tau\psi)_{(x,\tau,\eta)} = 
(d_\eta\psi)_{(x,\tau,\eta)} = 0  \},
$$
intersects the the space $\eta_1 = \ldots = \eta_n = 0$ transversally.
\item The map $(x,\tau,\eta) \mapsto (\frac{\del\psi}{\del \tau}, 
\frac{\del\psi}{\del \eta_1}
\ldots,  \frac{\del\psi}{\del \eta_{n}})$ has rank $n+1$ at every 
point of $C_\psi$, i.e. $\psi$ is non-degenerate.
\end{enumerate}

Define the map $F:C_\psi \to T^*M$ by $(x, \tau, \eta) \mapsto (x, 
(d_x\psi)_{(x,\tau,\eta)})$.  The image under $F$ of the subspace $\{\eta=0\} 
\subset C_\psi$ is homogeneous isotropic submanifold of
$T^*M$ of dimension $n+1$.  If this image is equal to $\Sigma$
then we say that $\psi$ \textit{parametrizes} $\Sigma$.
\begin{definition}
The space $I^m(M, \Sigma)$ of \textit{Hermite distributions} consists of 
consists of distributions on $M$ which have a local representation as
oscillatory integrals of the form
$$
\int e^{i\psi(x,\tau,\eta)} a\Bigl(x,\tau,\frac{\eta}{\sqrt\tau}\Bigr) d\tau 
\> d\eta , 
$$
where $\psi$ parametrizes $\Sigma$, and the amplitude $a(x,\tau, u)$  has the 
following properties (see 
\S3 of \cite{BG} for the precise formulation of the estimates):
\begin{enumerate}
\item $a(x,\tau, u)$ is rapidly decreasing as a function of $u$.
\item $a(x,\tau,u)$ is cutoff to be zero near $\tau=0$.
\item For sufficiently large $\tau$, $a(x,\tau,u)$ admits an expansion of the 
form
$$
a(x,\tau,u) \sim \sum_{i = 0}^\infty \tau^{m_i} a_i(x,u), 
$$
where each $m_i$ is either integer or half-integer, with $m_0 = m-1/2$ 
and $m_i \to -\infty$.
\end{enumerate}
\end{definition}

The conditions for a phase function to parametrize an isotropic ensure that a 
distribution in $I^m(M, \Sigma)$ will have wave-front set contained in $\Sigma$.

\subsection{The operator $A$}
We will now construct the operator $A$ that suitably approximates the operator
(\ref{0b}).  The idea is to construct a suitable parametrix of $D_\theta$ that
we will denote by $Q$.  

We start with the rolled-up version of the $\{\Box_k\}$, namely
\[
S = \Delta_H - (\trk)D_\theta.
\]
Although $S$ is not elliptic, the operator
$\Delta_P - (\trk)D_\theta = S+D_\theta^2$ is, and it commutes
with $S$.  We can therefore find an orthonormal basis of $L^2(P)$
of joint eigenfunctions of $S$ and $D_\theta$:
\[
S(\psi_j^k) = \lambda_j^k\,\psi_j^k,\quad 
D_\theta(\psi_j^k) = k\,\psi_j^k,\quad \lambda^k_1\leq \lambda^k_2\leq \cdots .
\]
Notice that, for each $k$, the $\lambda_j^k$, $j=1,2,\ldots$, are the
eigenvalues of $\Box_k$.
Therefore, the spectral drift phenomenon of Lemma \ref{Drift} has the following
interpretation in terms of the $\lambda$'s:  Denote the dimension
of $\calH_k$ by $d_k$.   Then, for all sufficiently large $k$
\begin{equation}\label{eitheror}
-\epsilon <\lambda_1^k\leq \cdots\lambda_{d_k}^k<\epsilon\quad\mbox{and}\quad
Mk\leq \lambda_{d_k+1}^k\leq \lambda^k_{d_k+2}\leq\cdots.
\end{equation}

Since the $\lambda_j^k+k^2$ are the eigenvalues of an elliptic
operator with positive symbol, only finitely-many of them
can be negative.  We can therefore form an operator that we'll denote by
\[
F:=\sqrt{S+D_\theta^2},
\]
because we can arrange that $F$ be a function of $S+D_\theta^2$ and
be such that $F^2-(S+D_\theta^2)$ is finite-rank.  The operator
$F$ is a standard elliptic first-order $\Psi$DO commuting with
$S$ and with $D_\theta$.  

\begin{lemma}\label{Cutoffs}
There exist non negative cut-off functions, $f,\ g\in C^\infty(\bbR)$,
such that the operator
\[
Q:= f(D_\theta ^2F^{-2})D_\theta^{-1} + g(S F^{-2}) F^{-1}
\]
has the following properties:
\begin{enumerate}
\item $Q$ is a classical elliptic $\Psi$DO of order $(-1)$.
\item In a neighborhood of the submanifold
$\calZ\subset T^*P\setminus\{ 0\}$,
the symbol of $Q$ is $\sigma(D_\theta)^{-1}$.
\item If $q_j^k$ denote  the eigenvalues of $Q$
(so that $Q(\psi_j^k) = q_j^k\psi_j^k$), then there exists
$c>0$ such that for all $k$ sufficiently large:
\begin{enumerate}
\item $q_j^k = k^{-1}$ for all $j=1,2,\ldots d_k$.
\item For all $j>d_k$  $q_j^k \geq c/k$.
\end{enumerate}
\end{enumerate}
\end{lemma}
\begin{proof}
We begin by noticing that the arguments of $f$ and $g$ above are
classical zeroth-order $\Psi$DOs, and therefore for
any choice of smooth functions $f$ and $g$ the operators
$f(D_\theta ^2F^{-2})$ and $g(S F^{-2})$ are classical
$\Psi$DOs of order zero.  If we choose 
$f(x)$ identically equal to zero near $x=0$ then
$f(D_\theta ^2F^{-2}) D_\theta^{-1}$ is a classical
$\Psi$DO of order $(-1)$, which implies (1).  Moreover,
if we choose $f(x)$ identically equal to one near $x=1$ then
condition (2) above is satisfied as well.  

For any $f$ and $g$ the eigenvalues of $Q$ can be
written as
\[
q = k^{-1}\, f\bigl(\frac{k^2}{\lambda+k^2}\bigr) +
\frac{1}{\sqrt{\lambda +k^2}}\,g\bigl(\frac{\lambda}{\lambda+k^2}\bigr),
\]
where we have written $\lambda$ for $\lambda_j^k$, for simplicity.
Notice that the arguments of $f$ and $g$ above add up to one.
Since $f(x)$ is being chosen identically equal to one near $x=1$, 
$f\bigl(\frac{k^2}{\lambda+k^2}\bigr) = 1$ if $k$ is large and $\lambda$
bounded.  If we choose $g(x)$ to be zero near $x=0$, this,
together with the previous choices, implies condition (3a).
It remains to be shown that, in addition, one can choose $f$
and $g$ so as to ensure condition (3b).

If $j>d_k$, then $\lambda_j^k>Mk$.  We distinguish two sub-regimes
of this case.  First, if $Mk<\lambda\leq ak^2$ for some $a>0$,
then
\[
\frac{k^2}{\lambda + k^2}\geq a/2,
\]
and so we can arrange for $f\bigl(\frac{k^2}{\lambda+k^2}\bigr)$
to be bounded below.  On the other hand, if $\lambda > ak^2$, then
$\frac{\lambda}{\lambda+k^2} > a$.  Therefore, if
$f(x)<1/10$ (for example) implies $g(1-x) = 1$, then $q$
is bounded below by a constant times $1/k$ as well.
\end{proof}

We now define the operator $A$ in terms of the operator $Q$ by:
\begin{equation}\label{defA}
A:=Q\circ S = Q\circ(\Delta_H - (\trk)D_\theta).
\end{equation}
Its main properties are summarized by the following:
\begin{corollary}\label{Adef}  \mbox{}

\begin{enumerate}
\item $A$ is a first-order classical $\Psi$DO with double characteristics,
$\calZ$.  The symbol of $A$ is identical to 
$\sigma(\Delta_H)/\sigma(D_\theta)$ in a conic neighborhood of $\calZ$.
\item If $\alpha_j^k$ denote the eigenvalues of $A$ 
(so that $A(\psi_j^k) = \alpha_j^k\psi_j^k$), then there exists
$C>0$ such that for all $k$ sufficiently large:
\begin{enumerate}
\item $\alpha_j^k = \frac{1}{k}\lambda_j^k
\in (-\epsilon/k,\,\epsilon/k)$ for all $j=1,2,\ldots d_k$.
\item For all $j>d_k$,  $\alpha_j^k \geq C$.
\end{enumerate}
\end{enumerate}
\end{corollary}
Since $A$ is a function of $S$, we can characterize the $\calH_k$
as the span of the joint eigenfunctions of $A$ and $D_\theta$ whose
eigenvalues are $O(1/k)$ from zero, with $A$ first-order.

\subsection{Proof of Theorem \ref{Main}.}
The family of operators $\{e^{-itA}\}_t$ is a smooth family
of Fourier integral operators, \cite{Tr} Ch.\ VIII \S 8.
The Schwartz kernel of the family, $U(t,x,y)$ is a distribution
on $\bbR\times P\times P$ with wave-front set equal to
\[
\Gamma := \{\,(t,\tau;x,\xi; y,-\eta)\;;\; \tau = p(x,\xi)\ 
\text{and}\ (y,\eta) = f_t(x,\xi)\,\}
\]
where $p$ is the principal symbol of $A$ and $f_t$ is the Hamiltonian 
flow of this symbol.  It follows by an
elementary wave-front set calculation that the wave-front set
of $\hphi(A)$ is contained in $\calZtr$.
Indeed the operator from $C^\infty(\bbR\times P\times P)$ to
$P\times P$ given by
\[
v(t,x,y)\mapsto \int_{\bbR}\,v(t,x,y)\,\phi(t)\,dt
\]
is a Fourier integral operator with canonical relation
\[
\calC := \{\,(t,\tau=0; x,\xi;y,\eta\,;\,x,\xi; y,\eta)\,\},
\]
and $\calC\circ\Gamma = \calZtr$.

Therefore, by the calculus of wave-front sets, to compute the Schwartz kernel
of $\hphi(A)$ modulo smooth functions we can microlocalize $U(t,x,y)$ to
any neighborhood of the set 
$ \{\,(t,\tau; x,\xi; y,-\eta)\;;\;\tau = 0\} $.
Let $\calV$ be a conic neighborhood of this set, specifically:
\[
\calV = \{ \,(t,\tau; x,\xi; y,-\eta)\;;\;|\tau|
< \delta\sqrt{\norm{\xi}^2+\norm{\eta}^2}\},
\]
where $\norm{\xi}$ is the norm with respect to the Kaluza-Klein metric on $P$.
Notice that (with a different constant, $\delta$)
\[
\Gamma_\delta :=
\calV\cap\Gamma =
\{\,(t,\tau;x,\xi; y,-\eta)\;;\; \tau = p(x,\xi)< \delta\norm{\xi}\ 
\text{and}\ (y,\eta) = f_t(x,\xi)\,\}
\]
\begin{lemma}
Let $W\subset P\times P$ be a neighborhood of the diagonal, and let $\epsilon >0$.
Then there exists a small enough $\delta>0$ such that the projection
\[
\begin{array}{ccc}
\Gamma_\delta &\to & \bbR\times P\times P \\
(t,\tau;x,\xi; y,-\eta) & \mapsto & (t,x,y)
\end{array}
\]
takes its values in $[-\pi-\epsilon, \pi+\epsilon]\times W$.
\end{lemma}
\begin{proof}
Let $\kappa$ denote the symbol of $D_\theta$, and let
\[
\norm{\cdot}^2 = \kappa^2 + h(\cdot)^2
\]
(that is, $h$ denotes the norm of the horizontal component of a covector).
Then, $p=h^2/\kappa$, and therefore
\[
 p< \delta \norm{\cdot} \quad \Leftrightarrow \quad
h^2 < \delta\,\kappa\,\sqrt{h^2+\kappa^2}.
\]
The flow $\{f_t\}$ is homogeneous of degree zero, therefore, to analyze
its behavior we can restrict our attention to the set: $\{\kappa = 1\}$.
But 
\[
\Bigl(\, \kappa = 1\ \text{and}\ p< \delta \norm{\cdot}\,\Bigr)
 \quad \Rightarrow \quad
\Bigl(\,\kappa=1\ \text{and}\ h^2<\delta(h+1)\,\Bigr),
\]
and $h^2<\delta(h+1)$ implies that $h = o(\delta)$.
Therefore, if $(t,\tau;x,\xi; y,-\eta)\in \Gamma_\delta\cap\{\kappa = 1\}$
then $h(x,\xi) = o(\delta)$, and in particular $(x,\xi)$ is close to
$\Sigma\cap\{\kappa = 1\}$.  Now notice that both $h$ and its Hamilton
vector field vanish on $\Sigma$.  Therefore, if $\delta$ is small enough,
the trajectory of $(x,\xi)$ by $\{f_t\}$ remains close to $(x,\xi)$ even
for times that are of the order of $\pi$.
\end{proof}
 
To prove that $\hphi(A)$ is an Hermite operator, recall (see \cite{Shu})
that  one can write the Schwartz kernel of $e^{-itA}$ as an oscillatory
integral with H\"ormander's phase function, i.e.\ of the form
\begin{equation}\label{oscint1}
U(t,x,y) = \int\,e^{i[\psi(t,x,y,\eta)- tp(y,\eta)]}\, a(t,x,y,\eta)\,
d\eta
\end{equation}
where $p$ is the principal symbol of $P$ and $\psi$ is homogeneous of
degree one in $\eta$ and satisfies:
\begin{equation}\label{oscint2}
\begin{array}{c}
p(x,d_x\psi(x,y,\eta)) = p(y,\eta)\\
\psi(x,y,\eta) = 0\quad\text{if}\ \inner{x-y}{\eta}=0\\
d_x\psi(x,y,\eta) = \eta\quad\text{if}\ x=y.
\end{array}
\end{equation}
The phase function $\psi$ parametrizes the conormal to the diagonal
on $P\times P$, and therefore the variable $\eta$ corresponds with a covector
at $y$.  The amplitude, $a(t,x,y,\eta)$ is of the form
\[
a(t,x,y,\eta)\sim\sum_{j=0}^\infty a_j(t,x,y,\eta)
\]
where $a_j$ is positive homogeneous in $\eta$ of order $(-j)$.  These
amplitudes must satisfy the transport equations which is always possible
for $t$ sufficiently small. 

Let us write down the symbol $p$, in our case.  Given a covector
$\eta$, write it as $\eta = (\eta_H,\kappa)$ with $\eta_H$ horizontal
and $\kappa$ a variable dual to $\partial_\theta$.  Then
\[
p(y,\eta) = \frac{\norm{\eta_H}^2}{\kappa}.
\]

Multiplying by $\phi(t)$ and integrating with respect to $t$
we obtain that the Schwartz kernel of $A_\phi$ is given by
oscillatory integrals of the form:
\begin{equation}\label{main}
\int e^{i\psi(x,y,\eta)}
f\left(x,y,\eta,\frac{\norm{\eta_H}^2}{\kappa} \right) d\eta\,d\kappa,
\end{equation}
where
$$
f(x,y,\eta,\kappa) = \int e^{-it\kappa} \phi(t) a(t,x,y,\eta)\ dt. 
$$
Since $\phi$ is compactly supported, $f$ is rapidly decreasing in
the $\kappa$ variable.  Thus, given the polyhomogeneous expansion
of $a(t,x,y,\eta)$, we see that (\ref{main}) is the local form of
an element of $I^0(P\times P, \calZtr)$.  (Technically $\psi$
parametrizes the conormal to the diagonal, rather than $\calZtr$.
But because of the rapid decay in the $\eta_H$ variables, specializing
to a phase function parametrizing $\calZtr$ would introduce only
smooth error terms.)

\section{The symbol computation}\label{sym.sec}

The symbol of a Hermite Fourier distribution in $I^l(M,\Sigma)$
is a \textit{symplectic spinor}, which is a half-density along $\Sigma$ tensored
with a smooth vector in the metaplectic representation associated to the
symplectic normal to $\Sigma$.  Denoting $\Sigma_\rho:= T_\rho\Sigma$, we can
write the space of symplectic spinors at $\rho$ as:
$$
\text{Spin}(\Sigma_\rho) = \hffm(\Sigma_\rho) \otimes H_\infty(\Sigma_\rho^{\circ}
/\Sigma_\rho).
$$

From Theorem \ref{Main} we know that $\hphi(A) \in I^0(P\times P, 
\calZtr)$.  Let $\phi\in \calZ$ and denote the tangent spaces
$$
V_\rho = T_\rho(T^*P), \qquad Z_\rho := T_\rho\calZ.
$$
so that $Z_\rho\subset V_\rho$ is a symplectic subspace.  The tangent space of
$T^*(P\times P)$ is $W_\rho := V_\rho\times V_\rho$, where the second factor
carries the opposite symplectic form (following the usual convention for 
canonical relations of FIO's).
And $Z^\triangle_\rho\subset W_\rho$ is thus isotropic.  With this set-up it is a 
simple exercise to work out that 
$$
(Z^\triangle_\rho)^{\circ}/(Z^\triangle_\rho) \isom  Z^{\circ}_\rho \oplus  
Z^{\circ}_\rho
$$
The smooth vector part of the symbol of $\hphi(A)$ is thus an element of 
$H_\infty(Z^{\circ}_\rho)\otimes H_\infty(Z^{\circ}_\rho)$.
By the identification $Z^\triangle_\rho \isom Z_\rho$ and the symplectic 
structure on the latter, $Z^\triangle_\rho$ carries a natural half-form.  
Thus we have a canonical identification
\begin{equation}\label{symident}
\text{Spin}(Z^\triangle_\rho) \isom \End(H_\infty(Z^{\circ}_\rho)).
\end{equation}

One last observation is needed to state the main result.  The metaplectic
group of $Z^{\circ}_\rho$ acts on $H_\infty(Z^{\circ}_\rho)$, by definition.  
Thus an element of the symplectic algebra of $Z^{\circ}_\rho$ acts on 
$H_\infty(Z^{\circ}_\rho)$ via the infinitesimal representation $d\tau$.
(We will describe this action in more detail below.) 
Now the projection $T_\rho(T^*P)\to T_{\pi(\rho)}P$
naturally identifies $Z^{\circ}_\rho$ with $T_{\pi(\rho)}X$.  Hence
the bundle automorphism $K$ defined by (\ref{Kdef}) can be lifted to an element 
of the symplectic
algebra of $Z^{\circ}_\rho$.  Under a canonical identification $Z^{\circ}_\rho
\isom \bbR^{2n}_{(x,\eta)}$ such that $H_\infty(Z^{\circ}_\rho) \isom 
\calS(\bbR^n_\eta)$ and 
$$
K = \begin{pmatrix} 0&(\kappa)\\ (-\kappa)& 0\end{pmatrix},
$$
(where $(\kappa)$ denotes the $n\times n$ diagonal matrix with entries
$\kappa_1,\dots,\kappa_n$), we have
$$
d\tau(K) = - \sum_{j=1}^n \kappa_j \Bigl( \frac{\partial^2}{\partial \eta_j^2} 
+ \eta_j^2
\Bigr),
$$
i.e. the harmonic oscillator Hamiltonian.  For this operator we have the
standard orthonormal basis of eigenfunctions $\ket{\nu} 
\in H_\infty(Z^{\circ}_\rho)$, $\nu = (\nu_1,\dots,\nu_n)$ such that
$$
d\tau(K)\ket{\nu} = (\kappa\cdot\nu+\trk/2)\ket{\nu},
$$
where $\kappa\cdot\nu = \sum_j \kappa_j \nu_j$.

\begin{proposition}\label{symbol.prop}
Under the identification (\ref{symident}), 
$$
\sigma(\hphi(A)) =  \sum_\nu \hphi(2\kappa\cdot\nu) \;\ket{\nu}\bra{\nu}.
$$
\end{proposition}

\subsection{The transport equation}
To prove Proposition \ref{symbol.prop}, we will use test distributions
in the class $I^l(P,\calR)$, associated to an isotropic ray $\calR = 
\{(p, r \alpha_p); r>0\}\subset \calZ$, for some fixed $p\in P$.
It suffices to make our calculation in the tangent space to $T^*P$
at the point $\rho = (p, \alpha_p)$ (because symbols will be homogeneous
in the ray coordinate).  Let $R_\rho = T_\rho\calR$, which is an isotropic subspace
of $Z_\rho$.  We can identify the symplectic normal space
$$
R_\rho^\circ/R_\rho \isom E_\rho \oplus \Zperp,
$$
where $E_\rho$ is the symplectic normal of $R_\rho$ as a subspace of $Z_\rho$.
Thus 
\begin{equation}\label{spinR}
\text{Spin}(R_\rho) \isom \hffm (R_\rho)\otimes H_\infty (E_\rho)
\otimes H_\infty(\Zperp).
\end{equation}

Consider the action of the operator $A$ on $I^l(P,\calR)$.  
From \cite{Gu}, Theorem V, we obtain the following 
transport equation:
\begin{proposition}\label{transport}
Since $\sigma(A)$ vanishes to second order on $\calR$, given $u\in I^l(P,\calR)$
we have $Au \in I^l(P,\calR)$ and
$$
\sigma(Au) = d\tau(\Hess(\sigma(A))) . \sigma(u) + \sigma_{sub}(A)\sigma(u),
$$
where the Hessian of the symbol, $\Hess(\sigma(A))$, acts on the smooth vector
component of  $\sigma(u_t)$ via the infinitesimal metaplectic representation
$d\tau$.
\end{proposition}

To describe the symbol of $\sigma(Au)$ more explicitly, we introduce convenient
local coordinates $(x,\theta; \eta,\tau)$ for $T^*P$ such that $\rho = (0,0;0,1)$.
Here $x$ is assumed to be a geodesic normal coordinate for $X$, so that, 
if $\beta$ denotes the metric on $X$, then $\beta_{ij} = \delta_{ij} + O(x^2)$.
The connection form is given by
$$
\alpha = d\theta + \sum \alpha_j(x) dx_j,
$$
and we may assume, without loss of generality, that $\alpha_j(0) = 0$ and
$\del_i\alpha_j(0) = \frac12 K_{ij}$, where $K_{ij}$ is the matrix of $K$
in these coordinates (which is the same matrix as $\omega$ since $g$ is normal
at 0).

The Kaluza-Klein metric on $P$, in block matrix form, is given by
\begin{equation}\label{kkg}
g = \begin{pmatrix}1&\alpha^t\\ \alpha&\beta + \alpha \alpha^t \end{pmatrix},
\end{equation}
(considering $\alpha_j$ to be a column vector) so that
\begin{equation}
g^{-1} = \begin{pmatrix}1 + \alpha^t\beta^{-1}\alpha&-\alpha^t\beta^{-1}\\
-\beta^{-1}\alpha& \beta^{-1} \\
\end{pmatrix}.
\end{equation}
From the connection between the symbol of $A$ and the symbol of 
$\Delta_H = \Delta_P - D_\theta^2$ in Corollary \ref{Adef} and the form of $g^{-1}$
we see that
$$
\sigma(A)(x,\theta;\eta,\tau) =  \frac{1}{\tau} (\eta - \tau\alpha)^t \beta^{-1} 
(\eta-\tau\alpha).
$$
The subprincipal symbol of $\Delta_H$ being zero, we also have 
$\sigma_{sub}(A) = -\trk$.
We need the Hessian of $\sigma(A)$, which is
\begin{equation*}
\begin{split}
\frac{\partial^2}{\partial x_i \partial x_j} \sigma(A)|_\rho &= 
2\sum \partial_i\alpha_k(0) \partial_j\alpha_k(0) = 
-\frac12(K^2)_{ij} \\
\frac{\partial^2}{\partial x_i \partial \eta_j} \sigma(A)|_\rho &= 
-2 \partial_i\alpha_j(0) = K_{ij} \\
\frac{\partial^2}{\partial \eta_i \partial \eta_j} \sigma(A)|_\rho &= 
2\delta_{ij} .\\
\end{split}
\end{equation*}
Thus, as an element of the symplectic algebra of $R_\rho^\circ/R_\rho$,
$$
\Hess(\sigma(A)) = \begin{pmatrix}K& 2I\\ \frac12 K^2 & K \end{pmatrix}.
$$

In order to express the transport equation in terms of the decomposition
\ref{spinR},
we need to express $\Hess(\sigma(A))$ as an element of the symplectic algebra 
of $E_\rho \oplus \Zperp$.  If we introduce coordinates:
\begin{equation*}\begin{split}
R_\rho^\circ/R_\rho &= \{(v,0;\xi,0)\} \\
E_\rho  &= \{(v_1,0;-(K\cdot v_1)/2,0)\}\\
\Zperp &= \{(v_2,0;(K\cdot v_2)/2,0)\}
\end{split}\end{equation*}
Then the canonical transformation $T: E_\rho \oplus \Zperp \to 
R_\rho^\circ/R_\rho$ is can be read off immediately:
\begin{equation*}\begin{split}
v &= v_1 + v_2 \\
\xi &= - (K\cdot v_1)/2 + (K\cdot v_2)/2.
\end{split}\end{equation*}
Conjugating the Hessian gives
\begin{equation*}\begin{split}
T^{-1}\circ\Hess(\sigma(A))\circ T 
&=\begin{pmatrix}1/2 &-K^{-1}\\ 1/2 & K^{-1} \end{pmatrix}
\begin{pmatrix}K& 2\\ K^2/2 & K \end{pmatrix}
\begin{pmatrix}1&1\\ -K/2 & K/2 \end{pmatrix}\\
&= \begin{pmatrix} 0& 0\\ 0 & 2K \end{pmatrix}
\end{split}\end{equation*}
Thus, by Proposition \ref{transport},
$\sigma(Au)$ is given in terms of the decomposition (\ref{spinR}) by 
the operator $2d\tau(K) - \trk$ acting on the 
$H_\infty(\Zperp)$ component of $\sigma(u)$.  If we take
$u \in I^l(P,\calR)$ with symbol
\begin{equation}\label{sigu}
\sigma(u) = \mu\otimes a \otimes \ket{\nu} \in \hffm (R_\rho)\otimes H_\infty 
(E_\rho) \otimes H_\infty(\Zperp),
\end{equation}
where $\{\ket\nu\}$ is the harmonic oscillator basis for $d\tau(K)$ introduced 
above, then the transport equation reduces to 
\begin{equation}
\label{sigAu}
\sigma(Au) = (2\kappa\cdot\nu)\sigma(u).
\end{equation}

\subsection{Proof of Proposition \ref{symbol.prop}}
With the transport equation understood, the proof of Proposition \ref{symbol.prop}
is fairly straightforward.  Suppose $u\in I^l(P,\calR)$ with symbol as in
\ref{sigu}, and let
$$
u_t = e^{-itA}u.
$$
Then $D_t u_t = -Au_t$, and hence, by the transport equation worked out above,
$$
D_t \sigma(u_t) = -2\kappa\cdot\nu\sigma(u_t).
$$
It is thus clear that
$$
\sigma(u_t) = e^{-2it\kappa\cdot\nu} \sigma(u),
$$
whence
$$
\sigma(\hphi(A)u) = \hphi(2\kappa\cdot\nu) \sigma(u).
$$
The result follows then from the interpretation of the smooth vector part of 
$\sigma(\hphi(A))$ as an endomorphism in (\ref{symident}).

\section{Constructing the projector}

To make the connection between $\hphi(A)$ and $\Pi$ we'll use a strategy
suggested to us by Victor Guillemin.
We first choose $\phi\in C_0^\infty$ with the following properties:
\begin{enumerate}
\item $\hphi(0)=1$ and for all integers $\ell\not=0$
$\hphi(\ell) = 0$.
\item $\hphi(\xi) < 1/2$ for all $\xi >\epsilon$, where
$\epsilon$ is the constant of Lemma \ref{Drift}.
\end{enumerate}
Such functions exist; for example we can take $\phi$ to be the convolution
of the characteristic function of $[-1,1]$ and a suitable function
in $C_0^\infty$.  A consequence of the first condition is that
the symbol of $\hphi(A)$ is, at each point in $\Sigma$, the
rank-one projector onto the ground state of the corresponding
harmonic oscillator.  Thus $\hphi(A)$ is a projector at the
symbolic level, and we will now use a Neumann series argument
to obtain the desired projector.

The basic idea is to use the identity 
\begin{equation}\label{heavy}
\frac{x-1/2}{\sqrt{1+4(x^2-x)}} + 1/2 = \Theta(x-1/2),
\end{equation}
where $\Theta$ is the Heaviside function.  Notice that,
by the second property of $\phi$ above, $\Theta(\hphi(A)-\frac{1}{2}I)$
is the desired projector.  Therefore we want to replace
$x$ by $\hphi(A)$ in (\ref{heavy}). 

We begin with
$E = 4(\hphi(A)^2 - \hphi(A))$.  Since the symbol of
$\hphi(A)$ is a projector, then $E$ is an Hermite FIO of
order $(-1/2)$, and therefore compact. 
Let us now define $(I+E)^{-1/2}$ by a Neumann series.
Following the remarks above, it is then a simple matter to check that
$$
\Pi := (\hphi(A)-\frac{1}{2}I) (I+E)^{-1/2} + \frac{1}{2}I
$$
is the desired projector.
It remains to show that $\Pi$ is Hermite.  The issue here is that
there is no guarantee that the limit of the
Neumann series for $(I+E)^{-1/2}$ produces an Hermite FIO.
However, all partial sums of the series are Hermite.
Let $(1+x)^{-1/2} = \sum_{l=0}^\infty c_l x^l$.
Then, by a standard Borel summation argument, 
one can construct an Hermite FIO,
$B \in I^0(P\times P, \calZ)$ such that for each $k$
$$
B - \sum_{l=0}^k c_l E^l \in I^{-\frac{k+1}{2}}(P\times P, \calZ).
$$
Let us now set
$$
R = \Pi - (\hphi(A)-\frac{1}{2}) B - \frac{1}{2}I.
$$
We will now show that $R$ is smoothing, which will imply that $\Pi$
is an Hermite FIO.  A calculation shows that
\[
R = (\hphi(A)-\frac{1}{2}I)\Bigl[(\sum_{l=0}^\infty c_l E^l) -B \Bigr].
\]
For each positive integer $k$, let us write
the operator in brackets in the form
\[
\sum_{l=0}^\infty c_l E^l -B = \Bigl(\sum_{l=0}^k c_l E^l - B\Bigr)
+ E^{k+1}\sum_{l=0}^\infty c_{l+k+1}\, E^l.
\]
and note that both terms map $H^0\to H^{(k+1)/2}$ (the series above
being bounded in $L^2$). Since this holds for any $k$,
$R$ is smoothing and $\Pi$ is Hermite.  Moreover, by construction.
$$
\sigma(\Pi) = \sigma(\hphi(A)) = \ket{0}\bra{0}.
$$
\hfill QED.

\end{document}